\documentclass{amsart}
\newtheorem{lemma}{Lemma}
  \newtheorem{thm}{Theorem}
   
    \newtheorem{remark}{Remark}
   
   \newtheorem{defn}{Definition}
  
\begin{document}
\title{ Computation of antieigenvalues of bounded linear operators via centre of mass}

\author{Kallol Paul, Gopal Das and  Lokenath Debnath}

\begin{abstract}

We introduce the concept of $\theta$-antieigenvalue and $\theta$-antieigenvector of a bounded linear 
operator on complex Hilbert space. We study the relation between $\theta$-antieigenvalue and centre of 
mass of a bounded linear operator and compute  antieigenvalue using the relation. This follows the notion of symmetric antieigenvalues introduced by Hossein et al. in \cite{19}. We show that the concept of real antieigenvalue, imaginary antieigenvalue and symmetric antieigenvalue follows as a special case of $\theta$-antieigenvalue. We also show how the concept of total antieigenvalue is related to the $\theta$-antieigenvalue. In fact, we show that  all the concepts of antieigenvalues studied so far follows from the concept of $ \theta$- antieigenvalue. We illustrate  with example  how to calculate the 
$\theta$-antieigenvalue  for an operator acting on a finite dimensional Hilbert space. 
\end{abstract}
\maketitle
\section{Introduction}
The concept of angle of an operator was introduced by Gustafson [5-7] in 1967 while studying the 
problems in perturbation theory of semi-group generators. For a bounded linear operator $T$ on a complex
 Hilbert space $H$ with the norm $\|.\|$ and the inner product $\langle,\rangle$ the cosine of the angle of the 
operator $T$ is defined as $$\cos \Phi (T) = \inf_{Tf \neq 0} \frac{Re\langle Tf,f \rangle}{\|Tf\| \|f\|}.$$ 
The similar concept was also introduced by Krei$\breve{n}$ \cite{20} in 1969 which he called the deviation 
of $T$ and denoted by $dev. T$. The $\cos \Phi (T)$ has another interpretation as first antieigenvalue ( real antieigenvalue) 
$\mu_1(T)$ of T which was also introduced by Gustafson \cite{8} : 
$$ \mu_{1}(T) = \inf_{Tf \neq 0} \frac{Re\langle Tf,f \rangle}{\|Tf\| \|f\|}. $$ 
The vector f for which $\mu_{1}(T)$ is attained (if exists), is called antieigenvector.
 The higher antieigenvalues are defined as follows : 
 $$ \mu_{n}(T) = \inf_{Tf \neq 0} \left\{ \frac{Re\langle Tf,f \rangle}{\|Tf\| \|f\|},~~~~~~~ f\bot \{f^{(1)}, f^{(2)}, f^{(3)},........,f^{(n-1)}\} \right\}$$
where $f^{(k)}$ is the $k^{th}$ antieigenvector of $T$.\\
The first antieigenvalue can be interpreted as the cosine of the largest angle ( real )  through which any vector can be rotated by the action of $T$. The concept of antieigenvalues is studied by Gustafson[8-12], Gustafson and Rao[13-14], Gustafson and Seddighin\cite{15}, Das et al\cite{2} ,Paul\cite{22}, Paul and Das\cite{23}. \\
Likewise Gustafson also introduced the concept of imaginary antieigenvalue as \\
$$ \inf_{Tf \neq 0} \frac{Im\langle Tf,f \rangle}{\|Tf\| \|f\|}.$$ 
The total cosine of an operator $T$ is defined as $$ | \cos | T = \inf_{Tf \neq 0 } \frac{|\langle Tf,f \rangle|}{\|Tf\|\|f\|} $$ 
The concept of total antieigenvalues is  studied by Gustafson and Seddighin \cite{16}, Seddighin\cite{24},  Hossein et. al.\cite{18}.  
In \cite{19}, first author Paul et.al. introduced the concept of symmetric antieigenvalue and antieigenvector for
 an operator $T$ as follows:-

$$ \Phi_{T}(f) =\frac{Re \langle Tf,f \rangle + Im \langle Tf,f \rangle}{\sqrt{2}\|Tf\|\|f\|}~~~~~~Tf \neq 0 $$
and  $$ \cos \Phi_{S}(T) = \inf_{Tf \neq 0} \Phi_{T}(f). $$

The symmetric antieigenvalue is also denoted by $\mu_{S}$. The vector $f$ for which $\Phi_{T}(f)$ attains
 the minimum (if exists ) is called the symmetric antieigenvector of $T$. For a self-adjoint operator $T$
 with the eigenvalues $\lambda_{1} \geq \lambda_{2} \geq \lambda_{3} \geq .......\geq \lambda{_n}$, 
$\cos \Phi(T) = \frac{\sqrt{\lambda_{1}\lambda_{n}}}{\lambda_{1}+\lambda_{n}}$, whereas $\cos \Phi(iT) =0$
[according to definition of Gustafson] rather abruptly, although $iT$ has the eigenvalues $i\lambda_{1}$, 
$i\lambda_{2}$, $i\lambda_{3}$,.....,$i\lambda_{n}$. But  the symmetric antieigenvalues of both $T$ and 
$iT$ are same i.e $\mu_{S}(T) = \mu_{S}(iT)$ for a self adjoint operator T. The definition of the 
antieigenvalue involves only real part of numerical range $W(T)$ of an operator $T$, but the symmetric 
antieigenvalue depends upon both real and imaginary part of the numerical range $W(T)$ of $T$. Following our concept of symmetric antieigenvalue Gustafson and seddighin \cite{17} studied slant antieigenvalues and slant antieigenvectors of operators, in which they didnot show any explicit relation between slant antieigenvalue and total antieigenvalue.\\
In \cite{23}, Paul and Das proved Min-max equality of a bounded linear operator $T$ on a complex Hilbert space $H$ using the concept of orthogonality of bounded linear operators in the sense of James\cite{3} and studied the relation between centre of mass and antieigenvalues. The Min-max equality in operator trigonometry goes as follows:\\

\textbf{Min-max equality:} \textit{For a bounded linear opeartor $T$ on a complex Hilbert space $$\sup_{\|x\|=1} \inf_{\epsilon \in R} \|(\epsilon T - I)x\|^{2} = \inf_{ \epsilon \in R}\sup_{\|x\|=1}\|(\epsilon T - I)x\|^{2}$$
 } 
The Min-max equality in operator trigonometry  was obtained by Gustafson\cite{4} in 1968, Asplund and Pt$\acute{a}$k\cite{1} in 1971. In \cite{23} using the concept of orthogonality of operators in the sense of James \cite{3} we  proved that   \\

\textit{For a bounded linear opeartor $T$ on a complex Hilbert space $$\sup_{\|x\|=1} \inf_{\lambda \in C} \|(\lambda T - I)x\|^{2} = \inf_{ \lambda \in C}\sup_{\|x\|=1}\|(\lambda T - I)x\|^{2}$$
 } 
In \cite{23} we introduced the concept of real centre of mass and total centre of mass of a bounded linear operator and studied their relation with antieigenvalues. We here mention the definitions of the real centre of mass and total centre of mass for the sake of completeness of information.\\

\begin{defn}
For any two bounded linear opeartors $T$ and $A$, there exists  scalars $\epsilon_{0} \in R $ and 
$\lambda_{0} \in C$ such that $$ \|T - \epsilon_{0}A\| \leq \|T - \epsilon A\|~~~~~~~~~\forall \epsilon 
\in R $$ and $$ \|T - \lambda_{0}A\| \leq \|T - \lambda A\|~~~~~~~~~\forall \lambda \in C $$
Then  numbers $\epsilon_{0}$ and $\lambda_{0}$ are called real centre of mass of $T$ with respect to $A$ and total 
centre of mass of $T$ with respect to $A$ respectively.
\end{defn}
We here introduce the concept of $\theta$-antieigenvalue and explore the relation between 
$\theta$-antieigenvalue and centre of mass of an operator. The concept of antieigenvalue, imaginary antieigenvalue and symmetric antieigenvalue then follows as a special case of $\theta$-antieigenvalue. We also show how the concept of total antieigenvalue is related to the concept of $ \theta$- antieigenvalue. Finally we give  examples of 
matrices to calculate $\theta$-antieigenvalue.

\section{$\theta$-antieigenvalues}

\begin{defn}[$\theta$-antieigenvalue]
Let $T$ be a bounded linear operator on a complex Hilbert space $H$ and $ \theta \in R .$
 Define $$\mu_{\theta} (f) = \frac{\cos \theta ~ Re\langle Tf,f \rangle + \sin \theta ~ Im\langle Tf,f \rangle}{\|Tf\|\|f\|},~~~~~Tf \neq \theta $$ and $$ \mu_{\theta}(T) = \inf_{Tf \neq 0} \frac{\cos \theta Re\langle Tf,f \rangle + \sin \theta Im \langle Tf,f \rangle}{\|Tf\|\|f\|} $$
$ \mu_{\theta}(T) $ is called $\theta$-antieigenvalue of $T$ and the vectors $f$ for which $\mu_{\theta}(T)$ attains the infimum (if exsits) are called $\theta$-antieigenvectors of $T$.
\end{defn}

If $\theta = 0$ then we get first antieigenvalue ( real antieigenvalue ) as 
$$\mu_{0}(T) = \inf_{Tf\neq 0}\frac{Re\langle Tf,f \rangle}{\|Tf\|\|f\|} = \cos T .$$
If $\theta = \frac{\pi}{2}$ then we get imaginary antieigenvalue  as 
$$\mu_{\frac{\pi}{2}}(T) = \inf_{Tf\neq 0}\frac{Im\langle Tf,f \rangle}{\|Tf\|\|f\|} $$
and if $\theta = \frac{\pi}{4}$ then we get symmetric antieigenvalue
$$ \mu_{\frac{\pi}{4}}(T)=\inf_{Tf\neq 0}\frac{Re\langle Tf,f \rangle+Im\langle Tf,f \rangle}{\sqrt{2}\|Tf\|\|f\|}=\mu_{s}(T). $$

We have studied the antieigenvectors in \cite{2,18,19,22} using the concept of stationary vectors, the definition of which
is given below :
\begin{defn}
{\bf Stationary vector.}\\
Let $ \phi(f) $ be a functional of a unit vector f $ \in $ H. Then
$ \phi(f) $ is said to have a stationary value  at f if the
function $ w_{g}(t) $ of a real variable t, defined as
\[ w_{g}(t) = \phi \left( \frac{f+tg}{\|f+tg\|}\right) \]
has a stationary value at t=0 for any arbitrary but fixed vector g
$ \in $ H.
The vector f is then called a stationary vector.\\
\end{defn}

We write \[ \Phi(f) = \frac{\cos \theta ~ Re\langle Tf,f \rangle + \sin \theta ~ Im\langle Tf,f \rangle}{\|Tf\|\|f\|}~ ; ~f \in
H,~ Tf \neq 0.  \] and find the necessary and
sufficient condition for a unit vector f to be a
stationary vector of $ \Phi(f) $.\\
For this we define
\[ w_{g}(t)~=~ \frac{ \left( \cos \theta ~ Re \langle T(f+tg),(f+tg)\rangle + \sin \theta ~ Im \langle T(f+tg),(f+tg)\rangle \right)^{2} }{\|T(f+tg)\|^2\|f+tg\|^2} \]
where g is an arbitrary but fixed vector of H.\\
If f is a stationary vector then we have $ w_{g}'(0) = 0 $ and so
we get
\[ 2 \|Tf\|^2( \cos \theta Af + \sin \theta Bf ) - ( a \cos \theta + b \sin \theta ) ( T^*T f + \|Tf\|^2 f ) = 0,\]
where $ A = Re T$, $ B = Im T $, $ a = Re \langle Tf,f \rangle $, $ b = Im \langle Tf,f \rangle $.

This is the necessary and
sufficient condition for $ \Phi(f) $ to be
stationary at a vector f. \\
Thus we obtain the following theorem
\begin{thm}
Let f be a unit $ \theta $-antieigenvector of  a bounded linear operator T. Then f satisfies the following equation
\[  2 \|Tf\|^2( \cos \theta Af + \sin \theta Bf ) - ( a \cos \theta + b \sin \theta ) ( T^*T f + \|Tf\|^2 f ) = 0,\]
where $ A = Re T$, $ B = Im T $, $ a = Re \langle Tf,f \rangle $, $ b = Im \langle Tf,f \rangle $.
\end{thm}
\begin{remark}
Putting $ \theta = 0, ~\frac{\pi}{2}, ~and~\frac{\pi}{4} $  we get the characteristic equation for f to be a real antieigenvector, imaginary antieigenvector and symmetric antieigenvector respectively.
\end{remark}
As we have studied in \cite{2,18,19} we can find  $\theta$-antieigenvalues for  bounded selfadjoint operators and normal operators and obtain similar results. We here find the relation between centre of mass and $ \theta$-antieigenvalue via Min-max equality. We first state the following theorem the proof of which is given in \cite{23} using the concept of orthogonality of bounded linear operators.
\begin{thm}
 For  bounded linear operators $T$ and $\lambda I$, where $\lambda \in C $ and $|\lambda|=1$, on a Hilbert space $H$
 $$ \sup_{\|f\|=1}\inf_{\epsilon \in R } \|(\epsilon T - \lambda I)f\|^{2} = \inf_{\epsilon \in R} 
\sup_{\|f\|=1}\|(\epsilon T - \lambda I)f\|^{2}$$
\end{thm}

Now we prove the following theorem.

\begin{thm}
Let $\epsilon_{0}$ be a real centre of mass of $T$ with respect to the opeartor $\lambda I$ 
[where $\lambda=\cos \theta +i\sin\theta, \theta \in R$] and $T$ is a bounded linear operator on a complex Hilbert space $H$. Then
 $$\mu_{\theta}(T) = \lim_{n \rightarrow \infty}\frac{\cos \theta Re\langle Tf_{n},f_{n}\rangle +
\sin \theta Im \langle Tf_{n},f_{n}\rangle }{\|Tf_{n}\|\|f_{n}\|}$$ where $\{f_{n}\}$ is a sequence of unit vectors in 
complex Hilbert space $H$ such that $Re \langle (\lambda I - \epsilon_{0}T)f_{n},Tf_{n}\rangle \rightarrow 0$ and $\|(\lambda I - \epsilon_{0} T)f_{n}\| \rightarrow \| \lambda I - \epsilon_{0} T\|.$
\end{thm}

\noindent \textbf{Proof:}
Since $\epsilon_{0}$ is a real centre of mass of $T$ with respect to $\lambda I$, we obtain a sequence of 
unit vectors $\{f_{n}\}$ in $H$ such that $Re \langle (\lambda I - \epsilon_{0}T)f_{n},Tf_{n}\rangle \rightarrow 0$ and 
$\|(\lambda I - \epsilon_{0} T)f_{n}\| \rightarrow \| \lambda I - \epsilon_{0} T\|$.\\
Now \begin{eqnarray*}
\|\lambda I - \epsilon_{0} T\| &=& \inf_{ \epsilon \in R} \|\lambda I - \epsilon T\|  \\
    &=& \inf_{\epsilon \in R } \sup_{\|f\|=1} \|(\lambda I - \epsilon T)f\|
\end{eqnarray*}
Using the min-max theorem we obtain 
\begin{eqnarray*}
\|\lambda I - \epsilon_{0} T\|^{2} &=& \lim_{n \rightarrow \infty} \left\{ 1 - 
\left(\frac{\cos \theta Re \langle Tf_{n},f_{n} \rangle + \sin \theta Im \langle Tf_{n},f_{n}\rangle }{\|Tf_{n}\|}\right)^{2}\right\}   \\
 &\leq& \sup_{\|f\|=1} \left\{1-\left(\frac{\cos \theta Re \langle Tf,f \rangle +
\sin \theta Im\langle Tf,f \rangle}{\|Tf\|}\right)^{2}\right\} \\
&\leq& \sup_{\|f\|=1}\inf_{\epsilon \in R} \|(\lambda I - \epsilon T)f\|^{2}   \\
&=& \|\lambda I - \epsilon_{0} T\|
\end{eqnarray*}

Therefore

\begin{eqnarray*}
\inf_{\|f\|=1} \frac{\cos \theta Re \langle Tf,f \rangle +\sin \theta Im \langle Tf,f \rangle }{\|Tf\|}  
&=& \lim_{n \rightarrow \infty} \frac{\cos \theta Re \langle Tf_{n},f_{n} \rangle + 
\sin \theta Im \langle Tf_{n},f_{n} \rangle }{\|Tf_{n}\|} \\
&=& \mu_{\theta}(T).
\end{eqnarray*}
This completes the theorem. \\
 \section{Total antieigenvalue and $\theta$-antieigenvalue}
We here show the relation between total antieigenvalue and $\theta$-antieigenvalue. 
\begin{lemma}
Let  $ f \in H,~ Tf\neq 0 $ and $ \theta \in R $. Then 
\[ \sup_{\theta \in R} \mu_{\theta} (f) = \sup_{\theta \in R} \frac{\cos \theta ~ Re\langle Tf,f \rangle + \sin \theta ~ Im\langle Tf,f \rangle}{\|Tf\|\|f\|} = \frac{|(Tf,f)|}{\|Tf\|\|f\|} .\]
\end{lemma}
\textbf{Proof.} For a fixed f we can think of $ \mu_{\theta} (f) $ as a function from R to R. Let $ \Psi : R \longrightarrow R $ be defined as
\[ \Psi(\theta) = \frac{\cos \theta ~ Re\langle Tf,f \rangle + \sin \theta ~ Im\langle Tf,f \rangle}{\|Tf\|\|f\|}. \]
Then using elementary calculus we see that $ \Psi $ attains its maximum at $ \cos  \theta  = \frac{ Re \langle  Tf, f \rangle}{| \langle Tf,f \rangle |} $ and $$ \sup_{\theta \in R} \Psi(\theta) = \frac{|\langle Tf,f \rangle|}{\|Tf\|\|f\|}. $$\\

\begin{thm}
$ | \cos | T = \inf_{Tf \neq 0} \sup_{\theta \in R} \mu_{\theta} (f).$
\end{thm}
\textbf{Proof.} Follows from Lemma 1 as 
\[ \sup_{\theta \in R} \Psi(\theta) = \sup_{\theta \in R} \frac{\cos \theta ~ Re\langle Tf,f \rangle + \sin \theta ~ Im\langle Tf,f \rangle}{\|Tf\|\|f\|} = \frac{|\langle Tf,f \rangle|}{\|Tf\|\|f\|} \] and so 
\[ | \cos | T = \inf_{Tf \neq 0} \sup_{\theta \in R} \mu_{\theta} (f) = \inf_{Tf \neq 0 } \frac{|\langle Tf,f \rangle|}{\|Tf\|\|f\|}.  \]
\begin{lemma}
$ \inf_{\theta \in R} \inf_{\epsilon \in R} \| \epsilon e^{i \theta}T - I \| = \inf_{\lambda \in C } \| \lambda T - I \|.$
\end{lemma}
\textbf{Proof.} By a result of \cite{21} there exists $ \lambda_0 = \epsilon_0 e^{i \theta_0} \in C $ such that $ \| \lambda_0 T - I \| = \inf_{\lambda \in C } \| \lambda T - I \|$ and so 
\[ \| \epsilon_0 e^{i \theta_0} T - I \| = \| \lambda_0 T - I \| \leq \| \lambda T - I \| = \| \epsilon e^{i \theta} T - I \| ~ \forall \epsilon \in R ~ and ~\theta \in R .\]
For each $ \theta \in R $ there exists $ \epsilon(\theta) = \epsilon_{\theta} $ such that 
\[ \| \epsilon_{\theta} e^{i \theta} T - I \| \leq \| \epsilon e^{i \theta} T - I \|~~ \forall \epsilon.\]
For $ \theta_0$, $ \epsilon_{\theta_0} $ may or may not be equal to $ \epsilon_0$ but 
\[ \| \epsilon_0 e^{i \theta_0} T - I \| \leq \| \epsilon_{\theta_0} e^{i \theta_0} T - I \| \leq \| \epsilon e^{i \theta_0} T - I \|~~ \forall \epsilon. \]
We choose 
\begin{eqnarray*}
 \epsilon(\theta) & = &  \epsilon_{\theta} ~~ \theta \neq \theta_0 \\
                  & = & \epsilon_0 ~~ \theta = \theta_0.
\end{eqnarray*}

Then 
\begin{eqnarray*}
& & \inf_{\theta \in R} \inf_{\epsilon \in R } \| \epsilon e^{i \theta}T - I \| \\
& = & \inf_{\theta \in R} \| \epsilon(\theta) e^{i \theta} T - I \| \\
& = & \| \epsilon_0 e^{i \theta_0} T - I \| \\
& = & \inf_{\lambda \in C } \| \lambda T - I \|.
\end{eqnarray*}
We next show that
\begin{thm} 
\[  \inf_{Tf \neq 0} \sup_{\theta \in R} \mu_{\theta} (f) = \sup_{\theta \in R}  \inf_{Tf \neq 0} \mu_{\theta} (f).\]
\end{thm}
\textbf{Proof.}  We have
\[  \cos T = \sqrt{ 1 - \inf_{\epsilon \in R} \| \epsilon T - I \|^{2}} ~so ~that~   \inf_{\epsilon \in R} \| \epsilon T - I \|^{2} = 1 - \inf_{Tf \neq 0}\left(\frac{Re\langle Tf,f \rangle}{\|Tf\|\|f\|}\right)^{2}\]
and 
\[  |\cos| T = \sqrt{ 1 - \inf_{\lambda \in C} \| \lambda T - I \|^{2}} ~so ~that~   \inf_{\lambda \in C} \| \lambda T - I \|^{2} = 1 - \inf_{Tf \neq 0}\left(\frac{|\langle Tf,f \rangle|}{\|Tf\|\|f\|}\right)^{2}\]

Using Lemma 2 we have 
\begin{eqnarray*}
\inf_{\lambda \in C}\|\lambda T-I\|^{2}& = & \inf_{\theta \in R}\inf_{\epsilon \in R}\|\epsilon e^{-i\theta}T-I\|^{2}  \\
\Rightarrow \inf_{\lambda \in C}\|\lambda T-I\|^{2} & = & \inf_{\theta \in R}\left\{  1 - \inf_{Tf \neq 0}\left(\frac{Re\langle e^{-i\theta}Tf,f \rangle}{\|e^{-i \theta}Tf\|\|f\|}\right)^{2} \right\} \\ 
\Rightarrow \inf_{\lambda \in C}\|\lambda T-I\|^{2} & = & 1 - \sup_{\theta \in R}\inf_{Tf \neq 0 } \mu_{\theta}^{2} (f) \\
\Rightarrow \sup_{\theta \in R}\inf_{Tf \neq 0 } \mu_{\theta}^{2} (f) & = &  1- \inf_{\lambda \in C}\|\lambda T-I\|^{2}\\
\Rightarrow \sup_{\theta \in R}\inf_{Tf \neq 0 } \mu_{\theta}^{2} (f) & = &  | \cos |^{2} T \\
\Rightarrow \sup_{\theta \in R}\inf_{Tf \neq 0 } \mu_{\theta} (f) & = & \inf_{Tf \neq 0} \sup_{\theta \in R} \mu_{\theta} (f).
\end{eqnarray*}
This completes the proof. \\

\noindent \textbf{Example } \\
$T=\left(
  \begin{array}{cc}
    2-3i & 0 \\
    0 & 3+2i \\
  \end{array}
\right)
$ be an operator on a two dimensional complex Hilbert space H. Let $z=(z_{1}, z_{2})^{t} \in H$, 
where $\left|z_{1}\right|^{2} + \left| z_{2}\right|^{2} = 1$. 
Then, $Re (Tz,e^{i\theta}z)=3\cos \theta + 2\sin \theta - (\cos \theta +5 \sin \theta)|z_{1}|^{2}$ 
 $\left|(Tz,z)\right|=\sqrt{\frac{26}{4}+26(|z_{1}|^{2} - \frac{1}{2})^{2}}$ and $\|Tz\|=\sqrt{ 13}$.
Then $|\cos | T =\frac{1}{\sqrt{2}}$ and the total cosine attains at the unit vector $z=(z_{1},z_{2})^{t}$ 
where $|z_{1}|^{2}= \frac{1}{2}$ and the total centre of mass is $\frac{5+i}{26}$.\\
The $\theta$-antieigenvalue of $T$ is as follows:-\\
\textbf{Case I:}  When $\sin(\theta + \tan^{-1}\frac{1}{5})\leq 0$, then $\mu_{\theta}(T) = \frac{3\cos \theta+2\sin\theta}{\sqrt{13}}$
 and it attains at the vector $z=(z_{1},z_{2})^{t}$ where $|z_{1}|=0$ and the centre of mass of $T$ with respect to $e^{i\theta} I$
where $I$ is the identity operator, is $\frac{3\cos \theta+2\sin\theta}{13}$\\
\noindent \textbf{Case II:} When $\sin(\theta + \tan^{-1}\frac{1}{5})>0$, then $\mu_{\theta}(T) = \frac{2\cos \theta-3\sin\theta}{\sqrt{13}}$
 and it attains at the vector $z=(z_{1},z_{2})^{t}$ where $|z_{1}|=1$ and the centre of mass of $T$ with respect to $e^{i\theta} I$
where $I$ is the identity operator, is $\frac{2\cos \theta-3\sin\theta}{13}$ \\
From $\theta$-antieigenvalue, we obtain antieigenvalue and symmetric antieigenvalue of $T$. When $\theta = 0$, then antieigenvalue $\cos T = 
\mu_{0}(T) = \frac{2}{\sqrt{13}}$ and the corresponding centre of mass is $\frac{2}{13}$. When $\theta=\frac{\pi}{4}
$ then the symmetric antieigenvalue is $\mu_{S}(T) = \mu_{\frac{\pi}{4}}(T) = -\frac{1}{\sqrt{26}}$ and the corresponding 
centre of mass $\epsilon_{0}=\frac{-1}{13\sqrt{2}}.$\\
Also $ \inf_{Tf \neq 0} \sup_{\theta \in R} \mu_{\theta} (f) = \sup_{\theta \in R}  \inf_{Tf \neq 0} \mu_{\theta} (f)
= \frac{1}{\sqrt{2}}. $ \\

\noindent {\bf ACKNOWLEDGEMENTS.}  We  thank   Professor T. K. Mukherjee for their invaluable
suggestions while preparing this paper.

\bibliographystyle{amsplain}

\begin{thebibliography}{99}

 \bibitem{1} Asplund E. and Pt$\acute{a}$k V., \textit{A Minimax Inequality for Operators and a Related Numerical Range}, Acta. Math, \textbf{126} (1971), 53-62.
 
\bibitem{2} Das K.C., DasGupta M. and Paul K., \textit{Structure of the antieigenvectors of a strictly accretive operator} International J. Math. and Math. Sci., \textbf{21} No. 4, (1998),761-766.

%\bibitem{3} Das K.C., Majumder S. and Sims B., \textit{Restricted numerical range and weak convergence on the boundary of numerical range}, Jour. Math. Phy. Sci., \textbf{21} No.1, (1987) , 35-42.

\bibitem{3} James R.C., \textit{Orthogonality and linear functionals in normed linear spaces}, Trans. Amer. Math. Soc. \textbf{61} (1947), 265-292.

\bibitem{4} Gustafson K., \textit{A Min-Max Theorem}, Notices Amer. Math. Soc. \textbf{15} (1968d), 799.

\bibitem{5} Gustafson K., \textit{Angle of an Operator and Positive Operator Products}, Bull. Amer. Math. Soc., \textbf{74} (1968a), 488-492.

\bibitem{6} Gustafson K., \textit{Positive(noncommutating) Operator Products and Semigroups}, Math. Zeit. \textbf{105} (1968b), 160-172.

\bibitem{7} Gustafson K., \textit{A note on left multiplication of semigroup generators}, Pacific J. Math. \textbf{24} (1968c), 463-465.

\bibitem{8} Gustafson K., \textit{Antieigenvalue Inequalities in Operator Theory }, Inequalities III, Proceedings Los Angeles Symposium, 1969 ed. O. Shisha, Academic Press  (1972), 115-119.

\bibitem{9} Gustafson K., \textit{ An extended operator trigonometry}, Linear Algebra Appl. \textbf{319} (2000), 117-135.
\bibitem{10} Gustafson K., \textit{ Operator Trigonometry}, Linear and Multilinear Algebra, \textbf{37} (1994), 139-159.

\bibitem{11} Gustafson K., \textit{Matrix Trigonometry}, Linear  Algebra Appl., \textbf{217} (1995), 117-140.

\bibitem{12}  Gustafson, K., \textit{Interaction antieigenvalues}, J. Math. Anal. Appl., \textbf{299} (2004),  174-185.


\bibitem{13}Gustafson K. and Rao D., \textit{Numerical Range and Accretivity of Operator Products}, J. Math. Anal. Appl., \textbf{60} (1977), 693-702.

\bibitem{14}Gustafson K. and Rao D., \textit{Numerical Range: The Field Values of Linear Operators and Matrices}, Springer, New York, 1997.

\bibitem{15} Gustafson K. and Seddighin M., \textit{Antieigenvalue Bounds}, J. Math. Anal. Appl., \textbf{143} (1989), 327-340.

\bibitem{16} Gustafson K. and Seddighin M., \textit{A note on total antieigenvalues}, J. Math. Anal. Appl., \textbf{178} (1993), 603-611.

\bibitem{17} Gustafson K. and Seddighin M., \textit{ Slant antieigenvalues and slant antieigenvectors of
operators }, Linear Algebra and its Applications, \textbf{432} (2010), 1348-1362.

\bibitem{18} Hossein Sk.M., Das K.C., Debnath L. and Paul K., \textit{Bounds for total antieigenvalue of a normal operator} International J. Math. and Math. Sci., \textbf{70}  (2004), 3877-3884.

\bibitem{19} Hossein Sk.M., Paul K., Debnath L. and Das K.C., \textit{Symmetric Anti-eigenvalue and Symmetric Anti-eigenvector}  J. Math. Analysis and Applications., \textbf{345}  (2008), 771-776.

\bibitem{20} Krei$\breve{n}$ H., \textit{Angular Localization of the Spectrum of a Mulplicative integral in a Hilbert space}, Functional Anal. Appl., \textbf{3}(1969), 89-90.

\bibitem{21}Paul K., Hossein Sk.M. and Das K.C., \textit{Orthogonality on B(H,H) and Minimal-norm Operator}, Journal of Analysis and Applications,   \textbf{6} (2008), 169-178.

\bibitem{22} Paul K., \textit{Antieigenvectors of the Generalized Eigenvalue Problem and an Operator Inequality Complementary to Schwarz's Inequality}, Novi Sad J. Math., \textbf{38} No. 2 (2008) 25-31.

\bibitem{23} Paul K. and Das G., \textit{Cosine of angle and center of mass of an operator},   Mathematica Slovaca, \textbf{62} (2012), No. 1, 109-122.

\bibitem{24} Seddighin M., \textit{ Antieigenvalues and total antieigenvalues of normal operators}, J. Math. Anal. Appl., \textbf{274} (2002) no.1, 239-254.
%\bibitem{23}Stampfli J.G., \textit{The Norm of A Derivation}, Pacific Journal of Mathematics, \textbf{33(3)} (1970), 737-747.


\end{thebibliography}

\bigskip
\noindent
\parbox[t]{.36\textwidth}{
Kallol Paul \\
Department of Mathematics\\
Jadavpur University \\
Kolkata 700032\\
INDIA\\
kalloldada@gmail.com\\ }
\parbox[t]{.40\textwidth}
{
Gopal Das\\
Department of Mathematics\\ Jalpaiguri Govt. Engg. College\\
 West Bengal 735102\\
  India\\
  gopaldasju@gmail.com\\}
\parbox[t]{.40\textwidth}
{
Lokenath Debnath\\
Dept of Mathematics\\
University of Texas-Pan American\\
1201 West University Drive\\
Edinburg, TX 78539\\
fax: (956) 384-5091\\
Phone: (956) 381-3459\\
E-mail: debnathl@utpa.edu\\
}

\end{document}